\documentclass{conm-p-l}
\copyrightinfo{2007}{The authors}

\usepackage{amssymb,amsmath} 
\usepackage{latexsym}

\newtheorem{theorem}{Theorem}[section]
\newtheorem{lemma}[theorem]{Lemma}
\newtheorem{proposition}[theorem]{Proposition}

\theoremstyle{definition}
\newtheorem{definition}[theorem]{Definition}

\newtheorem*{condition}{Condition BH}			

\theoremstyle{remark}
\newtheorem{remark}[theorem]{Remark}

\numberwithin{equation}{section}

\newcommand\mU{\mathcal{U}}

\def\buildrel#1_#2^#3{\mathrel{\mathop{\kern 0pt#1}\limits_{#2}^{#3}}}

\newcommand{\End}{\mbox{$\mathtt{End}$}}

\newcommand{\Der}{\mbox{$\mathtt{Der}$}}

\newcommand{\Ad}{\mbox{$\mathtt{Ad}$}}
\newcommand{\ad}{\mbox{$\mathtt{ad}$}}
\newcommand{\der}{\mbox{$\mathtt{Der}$}}
\DeclareMathOperator{\Span}{Span}

\def\cyclic{\mathop{\kern0.9ex{{+}
\kern-2.2ex\raise-.28ex\hbox{\Large\hbox
{$\circlearrowright$}}}}\limits}

\renewcommand{\a}{{\mathfrak{a}}{}} 
\renewcommand{\b}{{\mathfrak{b}}{}} 
\renewcommand{\r}{{\mathfrak{r}}{}}

\renewcommand{\k}{{\mathfrak{k}}{}} 
\renewcommand{\u}{{\mathfrak{u}}{}} 
 
\renewcommand{\S}{\mathbb S} 
\newcommand{\C}{\mathbb C}

\newcommand{\A}{\mathbb A}

\newcommand{\R}{\mathbb R}

\newcommand{\g}{{\mathfrak{g}}{}} 
\newcommand{\p}{{\mathfrak{p}}{}} 
\newcommand{\q}{{\mathfrak{q}}{}} 
\newcommand{\n}{{\mathfrak{n}}{}} 
\newcommand{\s}{{\mathfrak{s}}{}}

\newcommand{\so}{\mathfrak{so}}
\newcommand{\h}{{\mathfrak{h}}{}} 
\newcommand{\z}{{\mathfrak{z}}{}}

\newcommand{\CS}{{\mathcal S}{}} 
\newcommand{\CF}{{\mathcal F}{}} 
\newcommand{\CA}{{\mathcal A}{}} 
 
\newcommand{\CH}{{\mathcal H}{}}

\newcommand{\CD}{{\mathcal D}{}} 
\newcommand{\CM}{{\mathcal M}{}} 
\newcommand{\CR}{{\mathcal R}{}} 
\newcommand{\CK}{{\mathcal K}{}} 
\newcommand{\CN}{{\mathcal N}{}} 
\newcommand{\CE}{{\mathcal E}{}}

\newcommand{\mJ}{\mathrm{J}}

\begin{document}

\title[Quantization of anti de Sitter and symmetric spaces]
{Quantized Anti de Sitter spaces and non-formal \\ 
deformation quantizations of symplectic symmetric spaces}
\author[P. Bieliavsky]{Pierre Bieliavsky}
\address{\noindent Universit\'e catholique de Louvain, d\'epartement de 
math\'ematiques,
\qquad \hfill
Chemin du Cyclotron 2, B-1348 Louvain-la-Neuve, Belgium}
\email{bieliavsky  \, claessens \, voglaire \, @math.ucl.ac.be}
\thanks{The first author thanks Victor Gayral for enlightening discussions on
noncommutative spectral triples. This research is partially supported by 
the IAP grant ``NOSY'' at UCLouvain.}

\author[L. Claessens]{Laurent Claessens}
\thanks{Laurent Claessens is a FRIA-fellow.}

\author[D. Sternheimer]{Daniel Sternheimer}
\address{\noindent Institut de Math\'ematiques de Bourgogne, 
Universit\'e de Bourgogne, \qquad \hfill
BP~47870, F-21078 Dijon Cedex, France}
\curraddr{Department of Mathematics, Keio University, \\
3-14-1 Hiyoshi, Kohoku-ku, Yokohama, 223-8522 Japan}
\email{Daniel.Sternheimer@u-bourgogne.fr}

\author[Y. Voglaire]{Yannick Voglaire}
\subjclass[2000]{81S10; 53D50, 58B34, 81V25, 83C57}

\keywords{Non commutative geometry, BTZ-spaces, Deformation quantization, 
Symplectic symmetric spaces}

\begin{abstract}
We realize quantized anti de Sitter space black holes, building Connes spectral 
triples, similar to those used  for quantized spheres but based
on Universal Deformation Quantization Formulas (UDF) obtained from 
an oscillatory integral kernel on an appropriate symplectic symmetric space. 
More precisely we first obtain a UDF for Lie subgroups acting on a symplectic 
symmetric space $M$ in a locally simply transitive manner.  Then, observing 
that a curvature contraction canonically relates anti de Sitter geometry 
to the geometry of symplectic symmetric spaces, we use that UDF to define 
what we call Dirac-isospectral noncommutative deformations of the spectral 
triples of locally anti de Sitter black holes.  
The study is motivated by physical and cosmological considerations.
\end{abstract}

\maketitle
\vspace*{-3pt}

\section{Introduction}

\subsection{Physical and cosmological motivations}

This paper, of independent interest in itself, can also be seen as a small 
part in a number of long haul programs developed by many in the past decades, 
with a variety of motivations. The references that follow are minimal and 
chosen mostly so as to be a convenient starting point for further reading, 
that includes the original articles quoted therein. 

An obvious fact (almost a century old) is that anti de Sitter (AdS) space-time 
can be obtained from usual Minkowski space-time, deforming it by allowing 
a (small) non-zero negative curvature. The Poincar\'e group symmetry of 
special relativity is then deformed (in the sense of \cite{Ge64}) to the 
AdS group $SO(2,3)$. In $n+1$ space-time dimensions ($n \geq 2$) the 
corresponding AdS$_n$ groups are $SO(2,n)$. Interestingly these are 
the conformal groups of flat (or AdS) $n$ space-times. The \textsl{deformation 
philosophy} \cite{Fl82} makes it then natural, in the spirit of deformation 
quantization \cite{DS02}, to deform these further \cite{St05,St07}, i.e. 
quantizing them, which many are doing for Minkowski space-time. 

These deformations have important consequences. Introducing a small 
negative curvature $\rho$ permits to consider \cite{AFFS} massless particles 
as composites of more fundamental objects, the Dirac \textsl{singletons}, 
so called because they are associated with unitary irreducible representations 
(UIR) of $SO(2,3)$, discovered by Dirac in 1963, so poor in states that 
the weight diagram fits on a single line. These have been called Di and Rac 
and are in fact massless UIR of the Poincar\'e group in one space dimension 
less, uniquely extendible to the corresponding conformal group $SO(2,3)$ 
(AdS$_4$/CFT$_3$ symmetry, a manifestation of 't Hooft's holography). 
That kinematical fact was made dynamical \cite{FF88} in a manner consistent 
with quantum electrodynamics (QED), the photons being considered as 2-Rac 
states and the creation and annihilation operators of the naturally confined 
Rac having unusual commutation relations (a kind of ``square root" of the 
canonical commutation relations for the photon). 

Later \cite{FFS99} this phenomenon has been linked with the (then very 
recently) observed oscillations of neutrinos (see below, neutrino mixing).  
Shortly afterwards, making use of flavor symmetry, Fr\o nsdal was able to 
modify the electroweak model \cite{Fr00}, obtaining initially massless  
leptons (see below) that are massified by Yukawa interaction with Higgs 
particles. (In this model, 5 pairs of Higgs are needed and it predicts the 
existence of new mesons, parallel to the $W$ and $Z$ of the $U(2)$-invariant 
electroweak theory, associated with a $U(2)$ flavor symmetry.)

Quantum groups can be viewed \cite{BGGS} as an avatar of deformation 
quantization when dealing with Hopf algebras. Of particular interest 
here are the quantized AdS groups \cite{FHT93,Sta98}, especially at 
even root of unity since they have some finite dimensional UIR, a fact 
generally associated with compact groups and groups of transformations 
of compact spaces. It is then tempting to consider quantized AdS spaces 
at even root of unity $q=e^{i\theta}$ as ``small black holes" in an 
ambient Minkowski space that can be obtained as a limit when 
$\rho q \, \to \, -0$.  Note that, following e.g. 't Hooft 
(see e.g. \cite{tH06} but his approach started around 1980) that 
some form of communication is possible with quantum black holes by 
interaction at their surface. 

At present, conventional wisdom has it that our universe is made up mostly 
of ``dark energy'' (74\% according to a recent Wilkinson Microwave
Anisotropy Probe, WMAP), then of ``dark matter" (22\% according to WMAP), 
and only 4\% of ``our" ordinary matter, which we can more directly observe.  
Dark matter is ``matter'', not directly observed and of unknown composition,
that does not emit nor reflect enough electromagnetic radiation to be detected 
directly, but whose presence can be inferred from gravitational effects on 
visible matter. According to the Standard Model, dark matter accounts for 
the vast majority of mass in the observable universe.  Dark energy is a 
hypothetical form of energy that permeates all of space.  It is currently 
the most popular method for explaining recent observations that the universe 
appears to be expanding at an accelerating rate, as well as accounting for 
a significant portion of the missing mass in the universe.

The Standard Model of particle physics is a model which incorporates three 
of the four known fundamental interactions between the elementary particles 
that make up all matter (the fourth one, weakest but long range, 
being gravity). It came after the electroweak theory that incorporated QED 
(electromagnetic interactions) associated with the photon and the so-called 
weak interactions, associated with the \textsl{leptons} that now exist in three 
generations (flavors): electron, muon and tau, and their neutrinos.  
The Standard Model, of phenomenological origin, encompasses also the 
so-called strong interactions, associated with (generally) heavier particles 
called baryons (the proton and neutron, and many more), now commonly 
assumed to be bound states of ``confined" quarks with gluons, in 
three ``colors".  It contains 19 free parameters, plus 10 more in 
extensions needed to account for the recently observed \textsl{neutrino mixing} 
phenomena, which require nonzero masses for the neutrinos that are 
traditionally massless in the Standard Model. 

Very recently Connes \cite{CCM} developed an effective unified theory based  
on noncommutative geometry (space-time being the product of a Riemannian 
compact spin 4-manifold and a finite noncommutative geometry) for the 
Standard Model with neutrino mixing, minimally coupled to gravity.  
It has 4 parameters less and predicts for the yet elusive Higgs particle 
(responsible for giving mass to initially massless leptons) a mass that is 
slightly different (it is at the upper end of the expected mass range) from 
what is usually predicted. See also \cite{Co06}, and \cite{Ba06} in a 
Lorentzian framework.         

Previously, in part aiming at a possible description of quantum gravity, 
but mainly in order to study nontrivial examples of noncommutative manifolds,
Fr\"ohlich (in a supersymmetric context), then Connes and coworkers
had studied quantum spheres in 3 and 4 dimensions \cite{FGR99,CL01,CoDV}.  
The basic tool there is a spectral triple introduced by Connes \cite{Co94}. 
In the present paper we are developing a similar approach, but for hyperbolic 
spheres and using integral universal deformation formulas in the deformation 
quantization approach.

The distant hope is that these quantized AdS spaces (at even root of unity) 
can be shown to be a kind of ``small" black holes at the edge of our Universe 
in accelerated expansion, from which matter would emerge, possibly created as 
(quantized) 2-singleton states emerging from them and massified by interaction 
with ambient dark energy (or dark matter), in a process similar to those of 
the creation of photons as 2-Rac states and of leptons from 2-singleton 
states, mentioned above.

As fringe benefits that might explain the acceleration of expansion of our 
Universe, and the problems of baryogenesis and leptogenesis 
(see e.g. \cite{Cl06,Sa07}). Physicists love symmetries and even more 
to break them (at least at our level). One of the riddles that physics 
has to face is that, while symmetry considerations suggest that there 
should be as much matter as antimatter, one observes a huge imbalance 
in our region of the Universe. In a seminal paper published in 1967 that 
went largely unnoticed for about 13 years but has now well over a thousand 
citations (we won't quote it here), Andrei Sakharov addressed that problem, 
now called baryogenesis. If and when a mechanism along the lines hinted at 
above can be developed for creating baryons and other particles, it could 
solve that riddle. 

Roughly speaking the idea is that there is no reason, except theological, 
why everything (whatever that means) would be created ``in the beginning",
or as conventional wisdom has it now, in a Big Bang. There could very well 
be ``stem cells" of the primordial singularity that would be spread out, 
like shrapnel, mostly at the edge of the Universe. Our proposal is that 
these could be described mathematically as quantized AdS black holes.    
We shall now concentrate our study on them.

\subsection{Mathematical introduction}

Roughly speaking, a universal deformation formula (briefly UDF) for a 
given symmetry $\mathcal{G}$ is a procedure that, for every, say, topological 
algebra $\A$ admitting the symmetry $\mathcal{G}$, produces a deformation    
$\A_\theta$ of $\A$ within the {\sl same} category of topological algebras.  
Such a UDF is called {\sl formal} when the category it applies to is that    
of formal power series in a formal parameter with coefficients in associative 
algebras. 

For instance, Drinfel'd twisting elements in elementary quantum group theory 
constitute examples of formal UDF's (see e.g. \cite{ShPr}). Other formal 
examples in the Hopf algebraic context have been given by Giaquinto 
and Zhang \cite{GiZh}. In \cite{Za}, Zagier produced a formal example 
from the theory of modular forms.  The latter has been  used and generalized 
by Connes and Moscovici in their work on codimension-one foliations \cite{CoMo}.

In \cite{Ri}, Rieffel proves that von Neumann's oscillatory integral 
formula \cite{VonN} for the composition of symbols in Weyl's operator 
calculus actually constitutes an example of a non-formal UDF for the 
actions of $\R^d$ on associative Fr\'echet algebras. The latter has been
extensively used for constructing large classes of examples of noncommutative
manifolds (in the framework of Connes' spectral triples \cite{Co94})       
via Dirac isospectral deformations\footnote{A deformation triple 
$(\mathcal{A}_\theta, \mathcal{H}_\theta, D_\theta)$ is said 
\textbf{isospectral} when $\mathcal{H}_\theta$ and $D_\theta$ 
are the same for all values of $\theta$.} 
of compact spin Riemannian manifolds \cite{CL01} (see also \cite{CoDV}).   
Some Lorentzian examples have been investigated in \cite{BDRS} and 
\cite{Sitarz}. Other very interesting related approaches can be found in 
\cite{HNW,Ga05} and references quoted therein.                          

Oscillatory integral UDF's for proper actions of non-Abelian Lie groups 
have been given in \cite{Bi,BiMas,BiMae,BBM2}. Several of them were obtained 
through geometrical considerations on solvable symplectic symmetric spaces. 
Nevertheless, the geometry underlying the one in \cite{BiMas} remained unclear.

In the present work we build on these works in the AdS context, with when  
needed reminders of their main features so as to remain largely self-contained. 
First we show that the latter geometry is that of a   
solvable symplectic symmetric space which can be viewed as a curvature 
deformation of the rank one non-compact Hermitian symmetric space.  
[It can also be viewed as a curvature deformation of the AdS space-time, 
as we shall see in the last section of the article.]                           

Next, we develop some generalities on UDF's for groups which act strictly 
transitively on a symplectic symmetric space.  We give some precise criteria. 
We end the section by providing new examples with exact symplectic forms     
such as UDF's for solvable one-dimensional extensions of Heisenberg groups, 
as well as examples with non-exact symplectic forms.              

In the last section we apply these developments to noncommutative    
Lorentzian geometry.  In anti de Sitter space AdS$_{n\geq3}$, every  
open orbit $\CM_o$ of the Iwasawa component $\CA\CN$ of $SO(2,n)$ is 
canonically endowed with a causal black hole structure \cite{ClDe} 
(generalizing the BTZ-construction in dimension $n=3$). We define the analog  
of a Dirac-isospectral noncommutative deformation for a triple built on $\CM_o$. 
The deformation is maximal in the sense that its underlying 
Poisson structure is symplectic on the open $\CA\CN$-orbit $\CM_o$. 
In particular, it does not come from an application of Rieffel's 
deformation machinery for isometric actions of Abelian Lie groups. 
Moreover, via the group action, the black hole structure is encoded 
in the deformed spectral triple, with no other additional geometrical data,
in contradistinction with the commutative level\footnote{An interesting      
challenge would be to analyze which operator algebraic notions 
attached to the triple are responsible for the singular causality. 
That is not investigated in the present article.}.

\section{Curvature deformations of rank one Hermitian symmetric spaces 
and their associated UDF's}\label{RT}

\subsection{Preliminary set up and reminder}\label{sec:prelim}    
In \cite{BiMas}, a formal UDF for the actions of the Iwasawa component 
${\mathcal R}_0:=\CA\CN$ of $SU(1,n)$ is given in oscillatory integral form. 
It has been observed in \cite{BBM2} that this type of UDF is actually 
non-formal for proper actions on topological spaces. The precise framework   
and statement are as follows.  The group ${\mathcal R}_0$ is a               
one dimensional extension of the Heisenberg group  $\CN_0:=H_n$. 
Through the natural identification ${\mathcal R}_0=SU(1,n)/U(n)$ induced 
by the Iwasawa decomposition of $SU(1,n)$, the group ${\mathcal R}_0$ is 
endowed with a (family of) left-invariant symplectic structure(s) $\omega$. 
Denoting by $\r_0:=\a_0\times\n_0$ its Lie algebra, the map
\begin{equation}
\label{DARBOUX}
\r_0\longrightarrow {\mathcal R}_0:(a,n)\mapsto\exp(a)\exp(n)
\end{equation}
turns out to be a global Darboux chart on $({\mathcal R}_0,\omega)$.  
Setting $\n_0=V\times\R.Z$ with table $[(x,z)\,,\,(x',z')]=\Omega_V(x,x')\,Z$, 
and  $\r_0=\{(a,x,z)\,|\,,a,z\in\R;x\in V\}$, one has 

\begin{theorem}\label{BIMAS}
For all non-zero $\theta\in\R$, there exists a Fr\'echet function space 
$\CE_\theta$,                                                           
$C^\infty_c({\mathcal R}_0)\subset\CE_\theta\subset C^\infty({\mathcal R}_0)$,   
such that, defining for all $u,v\in C^\infty_c({\mathcal R}_0)$                  
\begin{equation}  \label{PRODUCT}
\begin{split}
(u\star_\theta v)(&a_0,x_0,z_0)\\
 &:=\frac{1}{\theta^{\dim {\mathcal R}_0}}
\int_{ {\mathcal R}_0\times {\mathcal R}_0}
\cosh(2(a_1-a_2))\,[\cosh(a_2- a_0)\cosh(a_0-a_1)\,]^{\dim\CR_0-2}\\
 & \times\exp\Big( \frac{2i}{\theta}\Big\{ S_V\big(\cosh(a_1-a_2)x_0, 
\cosh(a_2-a_0)x_1, \cosh(a_0-a_1)x_2\big)\\
 &- \cyclic_{0,1,2}\sinh(2(a_0-a_1))z_2 \Big\} \Big)\\
 &\times u(a_1,x_1,z_1)\,v(a_2,x_2,z_2)\, da_1da_2dx_1dx_2dz_1dz_2\;;
\end{split}
\end{equation}
where $S_V(x_0,x_1,x_2):=\Omega_V(x_0,x_1)+\Omega_V(x_1,x_2)+\Omega_V(x_2,x_0)$ 
is the phase for the Weyl product on $C^\infty_c(V)$ and                            
$\cyclic_{0,1,2}$ stands for cyclic summation\footnote{In \cite{BiMas},          
the exponent $\dim\CR_0-2$ was forgotten in the expression of the amplitude 
of the oscillating kernel.}, one has: 

\begin{enumerate}
\item  $u\star_\theta v$ is smooth and the map                                    
$ C^\infty_c({\mathcal R}_0)\times C^\infty_c({\mathcal R}_0)
\to C^\infty({\mathcal R}_0)$
extends to an associative product on $\CE_\theta$.                                          
The pair $(\CE_\theta,\star_\theta)$ is a (pre-$C^\star$) Fr\'echet algebra.
\item In coordinates $(a,x,z)$ the group multiplication law reads
\[ 
L_{(a,x,z)}(a',x',z')=\left(
a+a',e^{-a'}x+x',e^{-2a'}z+z'+\frac{1}{2}\Omega_V(x,x')e^{-a'}
\right).
\]
The phase and amplitude occurring in formula (\ref{PRODUCT}) are both invariant 
under the left action $L:{\mathcal R}_0\times {\mathcal R}_0\to {\mathcal R}_0$.
\item Formula (\ref{PRODUCT}) admits a formal asymptotic expansion of the form:
\begin{equation*}
u\star_\theta v\sim \,uv\,+\,\frac{\theta}{2i}\{u,v\}\,+O(\theta^2)\;;    
\end{equation*}
where $\{\,,\,\}$ denotes the symplectic Poisson bracket on 
$C^\infty({\mathcal R}_0)$ associated with $\omega$.
The full series yields an associative formal star product on 
$({\mathcal R}_0,\omega)$ denoted by $\tilde{\star}_\theta$.
\end{enumerate}
\end{theorem}

The setting and (i-ii) may be found in \cite{BiMas}, while (iii) is a            
straightforward adaptation  to ${\mathcal R}_0$ of \cite{BBM2}.    

\subsection{Geometry underlying the product formula}                             

We start with preliminary material concerning the symmetric spaces.              

\subsubsection{Symmetric spaces}

A {\it symplectic symmetric space} \cite{Bith,BCG} is a triple 
$(\mathcal{M},\omega,s)$ where $\mathcal{M}$ is a connected smooth manifold, 
$\omega$ is a non-degenerate two-form on $\mathcal{M}$ and 
$s:\mathcal{M}\times \mathcal{M}\to \mathcal{M}:(x,y)\mapsto s_x(y)$ is a 
smooth map such that $\forall x \in \mathcal{M}$ the map 
$s_x: \mathcal{M}\to \mathcal{M}$ is an involutive diffeomorphism of 
$\mathcal{M}$ preserving $\omega$ and admitting $x$ as an isolated fixed 
point. Moreover, one requires that the identity 
$s_x\circ s_y\circ s_x=s_{s_x(y)}$ holds for all $x,y \in \mathcal{M}$ \cite{Lo}.          
In this situation, if $x\in \mathcal{M}$ and $X$, $Y$ and $Z$ are smooth 
tangent vector fields on $\mathcal{M}$,                                                                         
\begin{equation}\label{CONNECTION}                                     
\omega_x(\nabla_XY,Z)\;:=\;\frac{1}{2}\,X_x.\omega(Y+s_{x_\star}Y\,,\,Z)
\end{equation}
defines an affine connection $\nabla$ on $\mathcal{M}$, the unique affine  
connection on $\mathcal{M}$ which is invariant under the symmetries 
$\{s_x\}_{x\in \mathcal{M}}$. 
It is moreover torsion-free and such that $\nabla\omega=0$.  
In particular, the two-form $\omega$ is necessarily symplectic.                  

It then follows that the group $\mathcal{G}=\mathcal{G}(\mathcal{M},s)$  
generated by the compositions $\{s_x\circ s_y\}_{x,y\in \mathcal{M}}$ is a 
Lie group of transformations acting transitively on $\mathcal{M}$. The group 
$\mathcal{G}$ is called the  {\bf transvection} group of $\mathcal{M}$. 
Given a base point $o$ in $\mathcal{M}$, the conjugation by the symmetry $s_o$ 
defines an involutive automorphism $\tilde{\sigma}$ of $\mathcal{G}$. 
Its differential at the unit element $\sigma:=\tilde{\sigma}_{\star_e}$ 
induces a decomposition into $\pm1$-eigenspaces of the Lie algebra $\g$ of 
$\mathcal{G}$: $\g=\k\oplus\p$.  The subspace $\k$ of fixed vectors turns out 
to be a Lie subalgebra which acts faithfully on the subspace $\p$ of 
``anti-fixed" vectors ($\sigma x = -x$ for $x \in \p$). Moreover $[\p,\p]=\k$. 
A pair $(\g, \sigma)$ as above is called a {\sl transvection pair}. 
The subalgebra $\k$ corresponds to the Lie algebra of the stabilizer of $o$ 
in $\mathcal{G}$, while the vector space $\p$ is naturally identified with the 
tangent space $T_o(\mathcal{M})$ to $\mathcal{M}$ at point $o$.  In particular 
the symplectic form at $o$, $\omega_o$, induces a $\k$-invariant symplectic 
bilinear form on $\p$. Extending the latter by $0$ on $\k$ yields a 
Chevalley 2-cocycle $\Omega$ on $\g$ with respect to the trivial representation
of $\g$ on $\R$.  A triple $(\g,\sigma, \Omega)$ as above is called a 
{\sl symplectic transvection triple}. It it said to be {\sl exact} when there 
exists an element $\xi$ in $\g^\star$ such that $\delta\xi=\Omega$, where 
$\delta$ denotes the Chevalley coboundary operator. Up to coverings, 
the correspondence which associates a symplectic transvection triple to a 
symplectic symmetric space is bijective. More precisely {\sl there is an 
equivalence of categories between the category of connected simply connected 
symplectic symmetric spaces and that of symplectic transvection triples} -- the 
notion of morphism being the natural one in both cases. In the above setting, 
exactness corresponds to the fact that the transvection group acts on 
$(\mathcal{M},\omega)$ in a strongly Hamiltonian manner.  

The above considerations can be adapted in a natural manner to the Riemannian or  
pseudo-Riemannian setting\footnote{See \cite{CaPa} for an excellent reference.}, 
essentially by replacing \emph{mutatis-mutandis} the symplectic structure by 
a metric tensor. The canonical connection (cf. Formula (\ref{CONNECTION}) above) 
corresponds in this case to the Levi-Civita connection.

\subsubsection{Symmetric spaces of group type}
We observe that in coordinates $(a,x,z)$ the map
$\phi:{\mathcal R}_0\to {\mathcal R}_0:(a,x,z)\mapsto(-a,-x,-z)$
preserves the symplectic form $\omega$ (because the coordinates are Darboux 
coordinates), is involutive and admits the unit element $e=(0,0,0)$ as an 
isolated fixed point. It may be therefore called ``centered symmetry''
of the associative kernel (\ref{PRODUCT}). Since the kernel is left-invariant 
as well, such a centered symmetry                                                              
$s_g: {\mathcal R}_0\to {\mathcal R}_0:g'\mapsto s_g(g')\;:=
L_g\circ\phi\circ L_{g^{-1}}\,(g')\;$ is attached to every point $g$ 
in ${\mathcal R}_0$. It turns out that endowed with the above family of 
symmetries, the manifold ${\mathcal R}_0$ becomes a 
symplectic symmetric space. More precisely:  

\begin{proposition}
Let $s:{\mathcal R}_0\times {\mathcal R}_0\to {\mathcal R}_0$ be the map        
$(g,g') \mapsto s_g(g')$.  Then the triple $({\mathcal R}_0,\omega,s)$ is a     
symplectic symmetric space whose transvection group is solvable. The underlying
affine connection is the unique affine symplectic connection which is invariant
under the group of symmetries of the oscillatory kernel (\ref{PRODUCT}).
\end{proposition}

\begin{proof}
In coordinates $(a,x,z)$, the symmetry map is expressed as 
\begin{equation}
\begin{split}
s_{(a,x,z)}(a',x',z')=(&2a-a',2  \cosh(a-a')x-x',\\
   			&2  \cosh(2 (a -  a'))z + \Omega_V(x,x') \sinh(a - a')-z').
\end{split}
\end{equation}
One then verifies that it satisfies the defining identity: 
$s_{g}\circ s_{g'}\circ s_g=s_{s_g(g')}$.  Concerning the solvability of the         
transvection group, four-dimensional symplectic transvection triples        
have been classified in \cite{Bith,Bi2}. The one we are concerned with here 
is given by Table (1) ($\varepsilon=1$) in Proposition~2.3 of \cite{Bi2};       
denoting here $\k=\Span\{ u_1,u_2,u_3 \}$ and 
$\p= \Span\{ e_1,e_2,f_1,f_2 \}$, it writes:  
\[                             
\begin{split}
[u_2,u_3]=u_1;\quad[u_2,e_2]=e_1;\quad [u_2,f_1]=-f_2;\quad[u_3,f_2]=e_1;
\quad[u_1,f_1]=2e_1;\\
[e_1,f_1]=2u_1;\quad[e_2,f_1]=u_3;\quad[e_2,f_2]=u_1;\quad[f_1,f_2]=u_2.
\end{split}  
\]
In these notations the Lie algebra $\r_0=\a_0\oplus V \oplus\R Z$
of ${\mathcal R}_0$ (with Darboux chart~(\ref{DARBOUX})) is generated by 
$\a_0 = \R.(-f_1)$, $V=\Span\{u_2-f_2, u_3+e_2\}$, $Z=2(u_1+e_1)$. 
We now verify (by a short computation) that for all 
$a,z\in\R$ and $x\in V$, one has 
\[ 
\tilde{\sigma}(\exp(-af_1)\exp(x+2z(u_1+e_1))K= (\exp(af_1)\exp(-x-2z(u_1+e_1))K,
\]
where $K$ denotes the analytic subgroup associated with $\k$.
Thus $\tilde{\sigma}$ gives $s_e=\phi$ on ${\mathcal R}_0$.
The triple considered here is therefore precisely the triple that induces on 
${\mathcal R}_0$ the present symmetric space structure, hence  
the transvection group is solvable. 
The higher dimensional case is similar to the 4-dimensional one. 

Note that the above symmetric space structure on ${\mathcal R}_0$ is            
canonically associated with the data of the oscillatory kernel (\ref{PRODUCT}).
Indeed, coming from the stationary phase expansion of an oscillatory integral, 
the formal star product $\tilde{\star}_\theta$ mentioned in item (iii) of      
Theorem \ref{BIMAS} is {\sl natural}  
in the sense that for all positive integer $r$ the $r$-th cochain of the 
star product is a bidifferential operator of order $r$. To every such natural 
star product is uniquely attached a symplectic connection \cite{Li, GuRa}.
In our case the latter, being invariant under the symmetries, must coincide 
with the canonical connection associated with the symmetries 
$\{s_g\}_{g\in {\mathcal R}_0}$.
\end{proof}
The above considerations lead to the following definitions.     
\begin{definition}
Let ${\mathcal R}$ be a Lie group. A {\bf symmetric structure} on 
${\mathcal R}$ is a diffeomorphism 
$\phi:{\mathcal R}\to{\mathcal R}$
such that 
$\phi^2=\mbox{id}_{\mathcal R}$; $\phi(e)=e$; 
$\phi_{\star_e}=-\mbox{id}_{T_e({\mathcal R})}$; and
setting,  for all $g\in{\mathcal R}$,
$s_g: =L_g\circ\phi\circ L_{g^{-1}}$  
then, for all $g'$, we define
$s_{g}\circ s_{g'}\circ s_g=s_{s_g(g')}$.   
\end{definition}
 Or equivalently:
\begin{definition}
A (symplectic) symmetric space, or more generally a homogeneous space, 
$\mathcal{M}$ of dimension $m$ is {\bf locally of group type} if there exists a 
$m$-dimensional (symplectic) Lie subgroup ${\mathcal R}$ of its automorphism 
group which acts freely on one of its orbits in $\mathcal{M}$.  
One says that it is {\bf globally} of group type if it is locally and    
if ${\mathcal R}$ has only one orbit\footnote{ In this case, for every choice 
of a base point $o$ in $\mathcal{M}$, the map 
${\mathcal R}\to \mathcal{M}:g\mapsto g.o$ is a diffeomorphism.}.
\end{definition}
Lie groups are themselves examples of symmetric spaces (globally) of group type. 
In the symplectic situation, however, a symplectic symmetric Lie group must 
be Abelian \cite{Bith}. We will see in what follows other non-Abelian examples.

\subsubsection{Strategy}
Our strategy for constructing UDF's for certain Lie groups may now be easily 
described: starting with a symplectic symmetric space of group type admitting 
an invariant deformation quantization --- obtained by geometric considerations 
at the level of the symmetric space structure --- one deduces a UDF for every 
Lie subgroup ${\mathcal R}$ as above by either identifying $\mathcal{M}$ 
to ${\mathcal R}$,  or, in the formal case, by restricting the deformation 
quantization to an open ${\mathcal R}$-orbit. This type of strategy, 
in contexts other than symmetric spaces, has already been proposed within    
a formal framework, see for example \cite{Xu,ShPr,Hung,CoMo,GiZh}.  

\section{Construction of UDF's for symmetric spaces of group type} 
The first example presented in the preceding section as well as the 
{\sl elementary solvable exact triples} (briefly ``ESET'') of \cite{Bi} are 
special cases of the following situation.

\subsection{Weakly nilpotent solvable symmetric spaces}
In this section $(\g,\sigma)$ denotes a complex solvable involutive Lie 
algebra (in short, iLa) such that if $\g=\k\oplus\p$ is the decomposition   
into eigenspaces of $\sigma$, the action of $\k$ on $\p$ is                 
nilpotent\footnote{This condition is automatic for a transvection algebra, but 
not in general. Indeed, consider the $2$-dimensional non-Abelian (solvable) 
algebra $\g=\mbox{Span}\{k,p\}$ with table $[k,p]=p$.}. 
We denote by $\mbox{pr}_\p:\g\to\p$ the projection onto $\p$ parallel to $\k$.

\begin{definition}
A {\bf good Abelian subalgebra} (in short, gAs) of $\g$ is an Abelian  
subalgebra $\a$ of $\g$, contained in $\p$, supplementary to a $\sigma$-stable 
ideal $\b$ in $\g$ and such that the homomorphism $\rho:\a\to\der(\b)$ 
associated with the split extension $0\to\b\to\g\to\a\to0$ is injective.  
\end{definition}

\begin{lemma}
If $(\g,\sigma)$ is not flat (i.e. $[\p,\p]$ acts nontrivially on $\p$) then 
a gAs always exists.                                                    
\end{lemma}

\begin{proof}
Since $\k$ is nilpotent, $[\k,\p]\neq\p$. Moreover by non-flatness there 
exists $X\in\p\backslash[\k,\p]$ not central in $\g$. Hence for every choice 
of a subspace $V$ supplementary  to $\a:=\R X$ in $\p$ and containing 
$[\k,\p]$, one has that $\b:=\k\oplus V$ is an ideal of $\g$ supplementary 
to $\a$ and on which $X$ acts nontrivially.
\end{proof}

Note that the centralizer $\z_\b(\a)$ of $\a$ in $\b$ is stable by the 
involution $\sigma$; indeed, for all $a\in\a$ and $X\in\z_\b(\a)$, one 
has $[a,\sigma X]=-\sigma[a,X]=0$.  Moreover, the map $\rho:\a\to \End(\b)$ 
being injective, we may identify $\a$ with its image: $\a=\rho(\a)$.  
Let $\Sigma:\End(\b)\to \End(\b)$ be the conjugation with respect to 
the involution $\sigma|_{\b}\in GL(\b)$, i.e. $\Sigma = Ad( \sigma|_{\b})$.  
The automorphism $\Sigma$ is involutive and preserves the canonical Levi 
decomposition $\End(\b)={\mathcal Z}\oplus \mathfrak{sl}(\b)$, where 
${\mathcal Z}$ denotes the center of $\End(\b)$.  Writing the element 
$a=\rho(a)\in\a$ as $a=a_Z+a_0$ within this decomposition, one has: 
$\Sigma(a)=a_Z+\Sigma(a_0)=-a=-a_Z-a_0$, because the endomorphisms $a$ 
and $\sigma|_{\b}$ anticommute.  Hence $\Sigma(a_0)=-2a_Z-a_0$ and therefore 
$a_Z=0$.  So, $\a$ actually lies in the semisimple part $\mathfrak{sl}(\b)$.  
For any $x\in \mathfrak{sl}(\b)$, we denote by                 
$x=x^S+x^N,\quad x^S,x^N\in \mathfrak{sl}(\b)$,                   
its abstract Jordan-Chevalley decomposition.  Observe that, denoting by 
$\mathfrak{sl}(\b)=\mathfrak{sl}_+\oplus \mathfrak{sl}_-$ the decomposition    
in $(\pm 1)$-$\Sigma$- eigenspaces, one has that $\a\subset \mathfrak{sl}_-$.   
Also, the algebra $\a_N:=\{ a^N\}_{a\in\a}$ is an Abelian subalgebra in 
$\mathfrak{sl}_-$ commuting with $\a$. Setting $\a_S:=\{a^S\}_{a\in\a}$, we have  

\begin{definition}\label{WNSA}
A gAs is called {\bf weakly nilpotent} if         
$\z_\b(\a_S)\subset\z_\b(\a_N)$ and\footnote{The last condition is automatic 
when $\b$ is Abelian. Observe also that it is satisfied when $\rho(\a)$ is 
contained in a Levi factor of the derivation algebra $\der(\b)$.} 
$\a_N\subset\der(\b)$.     
\end{definition}

Let  
$\b^c=:\bigoplus_{\alpha\in\Phi}\b_\alpha$   
be the weight space decomposition with respect to the action of $\a_S$. 
Note that for all $\alpha$, one has $\a_N.\b_\alpha\subset\b_\alpha$. 
Moreover, for all $X_\alpha\in\b_\alpha$ and $a^S\in\a_S$, one has
\begin{eqnarray*}
\sigma(a^S.X_\alpha)=\alpha(a^S)\sigma(X_\alpha)
=\sigma a^S\sigma^{-1}\sigma 
X_\alpha=\Sigma(a^S).\sigma(X_\alpha)=-a^S.\sigma(X_\alpha).
\end{eqnarray*}
Therefore, $-\alpha\in\Phi$ and $\sigma\b_\alpha=\b_{-\alpha}$.
Note in particular that $\sigma\b_0=\b_0$.
\begin{definition}
An involutive Lie algebra $(\g,\sigma)$ is called {\bf weakly nilpotent}  
if there exists a sequence of subalgebras $\{\a_i\}_{0\leq i\leq r}$
of $\g$ such that
\begin{enumerate}
\item $\a_0$ is a weakly nilpotent gAs of $\g$ with associated supplementary 
ideal $\b^{(0)}$.
\item $\a_{i+1}$ is a weakly nilpotent gAs of 
$\z_{\b^{(i)}}(\a_i)\quad(0\leq i\leq r-1)$ where, for $i\geq1$, $\b^{(i)}$  
denotes the $\sigma$-stable ideal of $\z_{\b^{(i-1)}}(\a_{i-1})$ associated 
with $\a_i$.                                                          
\item $\z_{\b^{(r)}}(\a_r)$ is Abelian.
\end{enumerate}
\end{definition}

\begin{proposition}
Assume the iLa $(\g,\sigma)$ to be weakly nilpotent. Then there exists 
a (complex) subalgebra $\s$ of $\g$ such that the restriction 
$\mbox{pr}_\p|_{\s}:\s\to\p$ is a linear isomorphism.
\end{proposition}

\begin{proof}
Let $\a$ be a weakly nilpotent gAs of $\g$ and set 
$V_\alpha:=\b_\alpha\oplus\b_{-\alpha}$ for every $\alpha\in\Phi$. 
Note that $V_0=\b_0$ and that each subspace $V_\alpha$ is then 
$\sigma$-stable and one sets $V_\alpha=\k_\alpha\oplus\p_\alpha$ for the 
corresponding eigenspace decomposition.   
Choose a partition\footnote{Such a partition can be defined as follows.  
Let $\h$ be a Cartan subalgebra of $\mathfrak{sl}(\b)$ containing $\a_S$ and   
let $\Lambda\in\h^\star$ denote the set of weights of the representation 
of $\mathfrak{sl}(\b)$ on $\b$. Note that the restriction map                  
$\lambda\mapsto\lambda|_{\a_S}$ from $\Lambda$ to $\a_S^\star$ is 
surjective onto $\Phi$. Let $\h=\h_\R\oplus i\h_\R$ be a real decomposition 
such that the restriction of the Killing form to $\h_\R$ is positive definite.
Every weight in $\Lambda$ is then real valued when restricted 
to $\h_\R$ \cite{Knapp}. Now, any choice of a basis of $\h_\R$ defines a 
partial ordering on $\Lambda$ with the desired properties. 
To pass to the set of weights $\Phi$, consider the $\h_\R$-components, 
$\a^\R_S$, of $\a_S$ viewed as a vector subspace of  $\h$. 
The restriction map $\rho:\Phi\to\Phi|_{\a^\R_S}$ is then a bijection. 
Indeed, for $\lambda\in\ker(\rho)$, one has, by $\C$-linearity, 
$\lambda(a+ia')=\lambda(a)+i\lambda(a')=0\,\forall a,a'\in\a_S^\R$. 
Hence $\lambda=0$ as an element of $\Phi$.  Therefore, an order on $\Lambda$ 
induces on $\Phi$ an order having the same properties.}           
of $\Phi\backslash\{0\}$ as $\Phi\backslash\{0\}=:\Phi^+\cup\Phi^-$ with the 
properties that $-\Phi^+=\Phi^-$ and that if $\alpha,\beta\in\Phi^+$         
with $\alpha+\beta\in\Phi$ then $\alpha+\beta\in\Phi^+$.  
One has $\b=\oplus_{\alpha\in\Phi^+}V_\alpha\oplus\b_0$. We set 
$\b^+:=\oplus_{\alpha\in\Phi^+}\b_\alpha$ and $\p^+:=\mbox{pr}_\p(\b^+)$.  
It turns out that the restriction map $\mbox{pr}_\p|_{\b^+}:\b^+\to\p^+$ 
is a linear isomorphism. Indeed for $X\in\b_\k:=\b\cap\k$ and $a\in\a_S$     
one has $\sigma(a.X)=\sigma a\sigma\sigma X=\Sigma(a).X=-a.X$; hence 
$\a_S.\b_\k\subset\p$. Therefore, for all $X\in\b_\alpha\cap\k\quad\alpha\neq0$, 
one can find $a\in\a_S$ such that $a.X=X\in\p\cap\k$; thus $\b_\alpha\cap\k=0$ 
as soon as $\alpha\neq0$, yielding $\ker(\mbox{pr}_\p|_{\b^+})=0$. 
The last condition of Definition \ref{WNSA} implies that $\a_S$ acts by     
derivations on $\b$, hence the usual argument yields that  $\b^+$ is a 
subalgebra normalized by $ \a\oplus\b_0$. Moreover the first condition      
of the same definition implies $[\a,\b_0]=0$.  The proposition follows by   
induction. One starts applying the above considerations to $\a=\a_0$ and 
$\b=\b^{(0)}$.  This yields a subalgebra $\s_0$ supplementary to 
$\z_{\b^{(0)}}(\a_0)$ and normalized by it. One then sets 
$\g_1:=\z_{\b^{(0)}}(\a_0)$, considers a weakly nilpotent $\a_1$ in $\g_1$ and 
gets a subalgebra $\s_1$ that is now supplementary to and normalized by   
$\z_{\b^{(1)}}(\a_1)$.  Applying this procedure inductively, one gets a 
sequence a subalgebras $\s_0,...,\s_{r-1},\s_r:=\z_{\b_r}(\a_r)\cap\p$ 
such that $\s_{i+1}$ normalizes $\s_i$ and $\s_i\cap\s_j=0$,
and defines $\s:=\oplus_i\s_i$.
\end{proof}

\begin{remark}
In the exact symplectic case, one has $\p_0\perp\a$ since $\xi[\a,\b_0]=0$.  
\end{remark}

\subsection{Darboux charts and kernels: an extension lemma}

Let $(\s_i,\Omega_i)$ with $i=1,2$ be symplectic Lie algebras, and denote 
by $({\bf S}_i,\omega_i)$ the corresponding simply connected symplectic Lie 
groups ($\omega_i$ being left-invariant). Given a homomorphism              
$\rho:\s_1\to\Der(\s_2)\cap\s\p(\Omega_2)$, form the corresponding          
semi-direct product 
$\s:=\s_1\times_\rho\s_2$   
and consider the associated simply connected Lie group ${\bf S}$ endowed 
with the left-invariant symplectic structure $\omega$  defined by 
$\omega_e:=\Omega_1\oplus\Omega_2$ ($e$ denotes the unit element in ${\bf S}$).

\begin{lemma}
If $\phi_i:(\s_i,\Omega_i)\to({\bf S_i},\omega_i)\quad(i=1,2)$ are Darboux  
charts, the map      
$\phi:(\s,\Omega:=\Omega_1\oplus\Omega_2)\longrightarrow({\bf S},\omega):
(X_1,X_2)\mapsto\phi_2(X_2).\phi_1(X_1)$       
is Darboux.
\end{lemma}

\begin{proof}
For $X \in \mathfrak{s}_1$ and  $Y \in \mathfrak{s}_2$ one has:
\[
\begin{split}
{L_{\phi^{-1}_{\star_\phi}}}\phi_\star(X)=
L_{\phi^{-1}_{\star_\phi}}\left(R_{{\phi_1}_{\star_{\phi_2}}}(\phi_{2_{\star}}X)\right)
&=L_{{\phi_1^{-1}}_\star}L_{{\phi_2^{-1}}_\star}R_{{\phi_1}_\star}{\phi_2}_\star X\\
&=\Ad\phi_1^{-1}\left(L_{{\phi_2^{-1}}_{\star_{\phi_2}}}{\phi_2}_{\star} X\right) 
\end{split}
\]
while for $Y\in\s_1$ one has                   
\[
{L_{\phi^{-1}}}_{\star_\phi}\phi_\star(Y)=
{L_{\phi_1^{-1}}}_\star  {L_{\phi_2^{-1}}}_\star {L_{\phi_2}}_\star{\phi_1}_\star Y=
{L_{\phi_1^{-1}}}_\star {\phi_1}_\star Y.
\]
Hence
\[
\omega_\phi(X,Y)=\Omega\Big(\Ad\phi_1^{-1}\big(
L_{{\phi_2^{-1}}_{\star_{\phi_2}}}{\phi_2}_{\star} X\big)\, ,\, 
{L_{\phi_1^{-1}}}_\star {\phi_1}_\star Y\Big)=0,
\]
because the first (resp. the second) argument belongs to $\s_2$ 
(resp. to $\s_1$), and for $X,X'\in\s_2$,
\[ 
\begin{split}
\omega_\phi(X,X')
&= \Omega\Big(\Ad\phi_1^{-1}\big( L_{{\phi_2^{-1}}_{\star_{\phi_2}}}{\phi_2}_{\star} X\big)\, ,
\, \Ad\phi_1^{-1}\big( L_{{\phi_2^{-1}}_{\star_{\phi_2}}}{\phi_2}_{\star} X'\big)\Big)\\
&= {\Ad\phi_1^{-1}}^\star\Omega_2\Big( L_{{\phi_2^{-1}}_{\star_{\phi_2}}}
{\phi_2}_{\star} X\, ,\, L_{{\phi_2^{-1}}_{\star_{\phi_2}}}{\phi_2}_{\star} X'\Big)\\
&= \Omega_2\Big( L_{{\phi_2^{-1}}_{\star_{\phi_2}}}{\phi_2}_{\star} X\, ,
\, L_{{\phi_2^{-1}}_{\star_{\phi_2}}}{\phi_2}_{\star} X'\Big)\\
&={\omega_2}_{\star_{\phi_2}}(X,X') \, \,  =\Omega_2(X,X').  
\end{split}  
\]
A similar (and simpler) computation applies for two elements of $\s_1$.
\end{proof}

A direct computation shows
\begin{lemma}\label{EXT}
Let $K_i\in \mbox{Fun}(({\bf S}_i)^3) $ be a left-invariant three point kernel 
on ${\bf S}_i$ ($i=1,2$).  Assume $K_2\otimes1\in \mbox{Fun}(({\bf S})^3)$ 
is invariant under conjugation by elements of ${\bf S}_1$.  Then 
$K:=K_1\otimes K_2\in\mbox{Fun}(({\bf S})^3)$ is left-invariant (under ${\bf S}$).
\end{lemma}

In particular, given associative kernels satisfying the above hypotheses,      
their tensor product defines an associative invariant kernel on the semidirect 
product ${\bf S}$. We will call it {\bf extension product} of $K_2$ by $K_1$. 

\subsection{Examples}

\subsubsection{One dimensional split extensions of Heisenberg algebras}
The strategy is as follows. First we observe that (the connected simply 
connected group associated to) every such extension acts simply transitively 
on the symmetric space $\mathcal{R}_0$, and is therefore diffeomorphic to it. 
Then, we remark that this diffeomorphism can be chosen to be symplectic. 
The latter is almost a homomorphism, up to an action of a subgroup of the 
automorphism group of the kernel on $\mathcal{R}_0$ (this is done by 
embedding the extension in a subgroup of the automorphism group of the 
kernel containing $\mathcal{R}_0$). Finally we check that the pullback of 
that kernel by the diffeomorphism gives an invariant kernel on the extension. 

So let us start with a Heisenberg Lie algebra $\mathfrak{h}_n=V\oplus \R Z$, 
where $(V,\Omega_V)$ is a $2n$-dimensional symplectic vector space, $Z$ 
is central, and $[v,v']=\Omega_V(v,v')Z$. One easily sees that a split
extension (with $i\circ p = \mathrm{Id}$) $\mathfrak{r}_n$ of               
$\mathfrak{h}_n$ by $\mathfrak{a}=\R A$:
\[
0 \rightarrow \mathfrak{h}_n \mathop{\rightarrow} \mathfrak{r}_n 
\mathop{\rightleftarrows}_i^p \mathfrak{a} \rightarrow 0 , 
\]
must be such that, for all $v\in V$ and $z\in \R$, 
$[A,v+zZ]=X.v+\mu(v)Z+2dzZ$   
with $X\in \mathrm{End}(V)$, $\mu\in V^*$ and $d\in\R$ such that 
$X-d\mathrm{Id}\in\mathfrak{sp}(V,\Omega_V)$. As symplectic form on the 
associated group, we choose the left-invariant 2-form whose value at the 
identity is the Chevalley 2-coboundary $\Omega=-\delta Z^*=Z^*[\cdot,\cdot]$, 
which is a natural generalization of that on $\mathcal{R}_0$ in \cite{BiMas}. 
This 2-coboundary is nondegenerate if and only if $d\neq 0$, which            
we will assume from now on\footnote{Quantizing extensions with $d=0$ 
requires a non exact 2-form. Since our method needs an exact one,        
we shall have to apply it on central extensions of our algebras. 
An example of that procedure is given in the next section.},                 
and consider extensions with parameters $dX$, $d\mu$ and $2d$, which we 
will denote by $(dX,d\mu,2d)$.

We have the following symplectic Lie algebra isomorphisms:
\begin{enumerate}
\item $(dX, d\mu, 2d) \cong (X,\mu,2)$ for all $d\in\R_0$, through the map 
$L(a,v,z)=(da,v,z)$.
\item $(X, \mu, d) \cong (X,0,2)$ for all $\mu\in V^*$, through the map 
$L(a,v,z)=(a,v+au,z)$ with $i_u\Omega_V=\mu$.
\item $(X, 0, 2) \cong (X',0,2)$ through the map $L(a,v,z)=(a,M.v,z)$, if and 
only if $M\in Sp(V,\Omega_V)$ is such that $MXM^{-1}=X'$, i.e. if 
$X-\mathrm{Id}$ and $X'-\mathrm{Id}$ belong to the same adjoint orbit 
of $Sp(V,\Omega_V)$.
\end{enumerate}
Thus we concentrate on algebras of type $(X,0,2)$ from which we can recover   
the quantization of the others. 

Now let $\mathfrak{r}_0=(I,0,2)$, $\mathfrak{r}'=(X,0,2)$, and 
$\mathcal{R}_0$, $\mathcal{R}'$ the associated groups. Elements of these 
algebras will be denoted respectively by $aA+v+zZ$ and $aA'+v+zZ$. 
The difference between the actions of $A'$ and $A$ on $V$ is 
$\bar X:=X-\mathrm{Id}$ which lies in $\mathfrak{sp}(V,\Omega_V)$. Extending the 
action $\bar X$ on $\mathfrak{r}_0$ by $[\bar X,aA+v+zZ]:=\bar X.v$, we 
can therefore view $\mathfrak{r}_0$ and $\mathfrak{r}'$ as subalgebras of 
the semidirect product $\mathfrak{g}=\mathfrak{r}_0\times\mathfrak{s}$ of 
$\mathfrak{r}_0$ by 
$\mathfrak{s}=\mathrm{span}(\bar X)\subset \mathfrak{sp}(V,\Omega_V)$. 
At the level of the groups ($\mathcal{S}$ corresponding to $\mathfrak{s}$), 
on the one hand we can identify $\mathcal{R}_0$ with $\mathcal{G}/\mathcal{S}$ 
as manifolds, and on the other hand as a subgroup of $\mathcal{G}$, 
$\mathcal{R}'$ acts on the quotient. That action is simply transitive. 
Indeed, on $\mathcal{R}_0$ and $\mathcal{R}'$, we have global coordinate maps 
$I(aA+v+zZ)=\exp(aA)\exp(v+zZ)$ and $I'(aA'+v+zZ)=\exp(aA')\exp(v+zZ)$   
such that $I'(aA'+v+zZ) = \exp(a\bar X).I(aA+v+zZ)$,      
giving a decomposition $g=sr\in\mathcal{R}'$ with $s\in\mathcal{S}$ and 
$r\in\mathcal{R}_0$. Thus, acting on $e\mathcal{S}\in\mathcal{G}/\mathcal{S}$, 
such an element $g$ gives            
$g\cdot e\mathcal{S} = C_s(r)\mathcal{S} = I(aA+e^{a\bar X}v+zZ)\mathcal{S}$, 
where $C_s$ denotes the action by conjugation, $C_s(g)=sgs^{-1}$.
The map $\phi:aA'+v+zZ\mapsto aA+e^{a\bar X}v+zZ$ is a diffeomorphism from 
$\mathfrak{r}'$ to $\mathfrak{r}_0$, and the corresponding map   
$\Phi=I\circ\phi\circ {I'}^{-1}$, or $sr\mapsto C_s(r)$, is a diffeomorphism 
from $\mathcal{R}'$ to $\mathcal{R}_0$. 
We now observe
\begin{proposition}
\begin{enumerate}
\item For $g=sr\in\mathcal{R}'$, we have 
$\Phi\circ L_g = L_{\Phi(g)}\circ C_s\circ\Phi$.
\item The kernel $K$ on $\mathcal{R}_0$ is invariant under the conjugations 
by $\mathcal{S}$.
\item Denoting by $\omega$ and $\omega'$ the left-invariant 2-forms with 
value $-\delta Z^*$ at the identity on $\mathcal{R}_0$ and $\mathcal{R}'$ 
respectively, we have $\Phi^*\omega=\omega'$.
\end{enumerate}
\end{proposition}
\begin{proof} 
\begin{enumerate}
\item For all $g=sr$, $g'=s'r'\in\mathcal{R}'$, we have 
$$\Phi\circ L_g(g') = \Phi(ss'C_{{s'}^{-1}}(r)r') 
= C_{ss'}(C_{{s'}^{-1}}(r)r') =L_{\Phi(g)}\circ C_s \circ\Phi(g').$$ 
\item Recall that for $s=e^{a'\bar X}$, 
$C_s(I(aA+v+zZ))=I(aA+e^{a'\bar X}v+zZ)$ and that $\bar X\in Sp(V,\Omega_V)$. 
Therefore in the kernel (\ref{PRODUCT}) the amplitude is independent of $v$, 
and in its phase, the function $S_V$ is invariant under the symplectomorphisms  
of $V$. The kernel as a whole is thus also invariant under $Sp(V,\Omega_V)$.
\item By left-invariance, the condition is 
$\omega_e\circ{\left(L_{\Phi(g)^{-1}}\circ\Phi\circ L_g\right)_*}_e=\omega'_e$. 
Using the first property of $\Phi$ above and the invariance of $\omega_e$    
under $Sp(V,\Omega_V)$, we get   
$\omega_e\circ{\left(L_{\Phi(g)^{-1}}\circ\Phi\circ L_g\right)_*}_e 
=\omega_e\circ {\left(C_s\circ\Phi\right)_*}_e    
=\omega_e\circ {\phi_*}_0$,      
which is readily seen to be equal to $\omega'_e$.
\end{enumerate}
\end{proof}

Defining a kernel 
$K':\mathcal{R}'\times\mathcal{R}'\times\mathcal{R}'\rightarrow \C$ 
on $\mathcal{R}'$ by $K'=\Phi^*K$, we now have 
\begin{proposition} The kernel $K'$ is
\begin{enumerate}
\item invariant under the diagonal left action of $\mathcal{R}'$,
\item associative.
\end{enumerate}
Together with the functional space $\Phi^*\mathcal{E}_\theta$, it thus 
defines a WKB quantization of $\mathcal{R}'$.
\end{proposition}

A quantization of these groups can be obtained by the same method as 
in \cite{BiMas}. So let us first quickly review that method.  
On a connected simply connected symplectic solvable Lie group 
$(\mathcal{R}, \omega)$, with $\omega$ a left-invariant exact symplectic 
form (so that the action by left translations is strongly Hamiltonian), 
one chooses a global Darboux chart $I:\mathfrak{r}\rightarrow\mathcal{R}$ 
for which the Moyal star product is \emph{covariant}, i.e. if 
for $X\in\mathcal{r}$, $\lambda_X\in C^\infty(\mathfrak{r})$ denotes the 
(dual) moment map in these coordinates, and $\star^M_\theta$ the Moyal product 
on $\mathfrak{r}$, one has 
$[\lambda_X,\lambda_Y]_{*^M_\theta} = 2\theta\{\lambda_X,\lambda_Y\}$.  
Such charts always exist on these groups (see \cite{Ped}, \cite{ArCo}). 

Covariance of the Moyal product implies that 
$\rho_\theta:\mathfrak{r}\rightarrow \mathrm{End}(C^\infty(\mathfrak{r})
[[\theta]]):X\mapsto [\lambda_X,\cdot]_{*_\theta^M}$ 
is a representation of $\mathfrak{r}$ by derivations of the algebra 
$(C^\infty(\mathfrak{r})[[\theta]],\star_\theta^M)$. In order to find an 
invariant product on $\mathfrak{r}$, one then tries to find an invertible 
operator $T_\theta$ which intertwines this action and that by fundamental 
vector fields, i.e. such that 
$T_\theta^{-1}\circ \rho_\theta(X)\circ T_\theta = X^*$. 
Those found up to now were all integral operators of the type 
$\mathcal{F}^{-1}\circ\phi_\theta^*\circ\mathcal{F}$, where $\mathcal{F}$ 
is a partial Fourier transform and $\phi_\theta$ a diffeomorphism (see the 
proof of Theorem~\ref{thmL2} for the precise form of the one in \cite{BiMas}). 
An invariant product $\star_\theta$ on $\mathfrak{r}$ is then defined by 
$u\star_\theta v=T_\theta^{-1}(T_\theta u\star_\theta^W T_\theta v)$, where 
$\star_\theta^W$ is the Weyl product on $\mathfrak{r}$. Modulo some work on   
the function spaces, this gives a WKB invariant quantization of $\mathcal{R}$.
In our case, choosing as Darboux chart the map 
$I(aA+v+zZ)=\exp(aA)\exp(e^{-a\bar X}v+zZ)$, 
one checks that the same integral operator $T_\theta$ as in \cite{BiMas} 
works here, giving thus rise to the same kernel (\ref{PRODUCT}). 
This reflects again the fact that our groups are all subgroups of the 
automorphism group of a symplectic symmetric space on which they act 
simply transitively.

\subsubsection{Non exact example}

As mentioned before, all the examples of symplectic symmetric Lie groups shown 
up to now were endowed with an exact symplectic form, as our method requires 
exactness in order for the left translations to be strongly Hamiltonian. 
We present here an example with a non exact symplectic form showing that, 
as expected, considering a central extension allows to apply the same method.

Let $\mathfrak{r}=\mathrm{span}\left<A,V_1,V_2,W_1,W_2,Z\right>$ be the         
Lie algebra defined by                                     
$[A,V_i]=V_i$, $[A,W_i]=-W_i$, $[A,Z]=2Z$, $[V_1,V_2]=Z$,
and choose the non exact\footnote{One can actually show that every symplectic 
2-cocycle on $\mathfrak{r}$ is non exact and that what follows can be applied 
to any one of them.} Chevalley 2-cocycle $\Omega$:              
$\Omega(A,Z)=1$, $\Omega(V_i,W_i)=1$, $\Omega(V_1,V_2)=1/2$.         
We define $\mathfrak{g}$ as the central extension of $\mathfrak{r}$ by the 
element $E$ with commutators                                    
$[X,Y]_{\mathfrak{g}} = [X,Y]_{\mathfrak{r}} + \Omega(X,Y).E$  
for all $X,Y\in\mathfrak{g}$, where we extended $[\cdot,\cdot]_{\mathfrak{r}}$ 
and $\Omega$ by zero on $E$. Then $\Omega=-\delta E^*$ is a 2-coboundary.      

Now the connected simply connected Lie group $\mathcal{R}$ whose Lie algebra 
is $\mathfrak{r}$ can be realized as the coadjoint orbit $\mathcal{O}$ of 
$E^*$ in $\mathfrak{g}^*$, and a global Darboux chart $J$ from $\R^6$     
to $\mathcal{O}$ is given by:                                             
\begin{eqnarray*}
J(q_1,p_1,q_2,p_2,q_3,p_3)&=&\exp(q_3H) \exp((p_3+q_1p_1+q_2p_2-p_1p_2/2)Z) \\
&&\exp((q_2+p_1/4)W_1) \exp(p_2V_1) \\
&&\exp((q1-p_2/4)W_2) \exp (p_1V_2) \cdot E^*.
\end{eqnarray*}
In this chart, the (dual) moment maps are linear in the $p_i$, so that the   
Moyal product is covariant. From now on, the method outlined above can be    
carried on the same way as in \cite{BiMas}, with the same integral operator 
$T_\theta$ as before.

\subsubsection{The Iwasawa factor of $\mathfrak{sp}(2,\R)\simeq \so(2,3)$}  
\label{SubSubSecIwasawaSP}

The Lie algebra $\g_0:=\mathfrak{sp}(n,\R)$ of the group $\mathcal{G}_0:=Sp(n,\R)$ 
is defined as the set of $2n\times 2n$ real matrices $X$ such that 
${}^\tau XF+FX=0$ where 
$F:=\left(\begin{array}{cc}0 &I_n\\-I_n&0\end{array}\right)$.  One has
\begin{equation}
\g_0=\{\left(\begin{array}{cc}A &S_1\\S_2&-{}^\tau A\end{array}\right) 
\mbox{ where }A\in\mbox{Mat}(n\times n,\R)\mbox{ and }
S_i={}^\tau S_i\,;i=1,2\}. 
\end{equation}
In particular, $F\in \mathcal{G}_0$ and a Cartan involution of $\g_0$ is given 
by $\theta:=\Ad(F)$.  The corresponding Cartan decomposition 
$\g_0=\k_0\oplus\p_0$ is then given by
\begin{equation}
\k_0\simeq\u(n)\mbox{ and }
\p_0=\{\left(\begin{array}{cc}S&S'\\S'&-S\end{array}\right)\},
\end{equation}
where the matrices $S$ and $S'$ are symmetric. For $n=2$, a maximal Abelian 
subalgebra $\a$ in $\p_0$ is generated by   
$H_1=E_{11}-E_{33}$ and $H_2=E_{22}-E_{44}$          
where as usual $E_{ij}$ denotes the matrix whose component are zero except   
the element $ij$ which is one. The restricted roots $\Phi$ w.r.t $\a$ are 
then given by 
\begin{equation}
\Phi=\{\alpha_0,\;\alpha_1,\;\alpha_2:=\alpha_0+\alpha_1
,\;\alpha_3:=\alpha_0+2\alpha_1\}
\end{equation}
 with $\alpha_0:=2H_2^*$, $\alpha_1:=H_1^*-H_2^*$, 
hence  $\alpha_2:=H_1^*+H_2^*$ and $\alpha_3:=2H_1^*$,
 where $H_i^*(H_j):=\delta_{ij}$. The corresponding root spaces $\g_{\alpha_i}$
$(i=0,...,3)$ are one-dimensional, generated respectively by 
$N_0=E_{24}$, $N_1=E_{12}-E_{43}$, $N_2=E_{14}+E_{23}$, $N_3=E_{13}$. 
With this choice of generators, the minimal parabolic subalgebra 
$\s:=\a\oplus\n$ with $\n:=\bigoplus_{i=0}^3\g_{\alpha_i}$ has the 
following multiplication table:
\begin{subequations}\label{EqTableSP}
\begin{align} 
[H_1,N_1]&=N_1,			&\left[H_1,N_2\right]&=N_2,\\
\left[H_2,N_1\right]&=-N_1,	&\left[H_2,N_2\right]&=N_2,\\
\left[H_1,N_3\right]&=2N_3,	&\left[H_2,N_0\right]&=2N_0,\\
\left[N_0,N_1\right]&=-N_2,	&\left[N_1,N_2\right]&=2N_3;
\end{align}
\end{subequations}
the other brackets being zero.  Setting 
\begin{equation}
\s_1:=\Span\{H_2,N_0\}\mbox{ and }\s_2:=\Span\{H_1,N_1,N_2,N_3\},
\end{equation}
one observes that $\s$ is a split extension of $\s_2$ by $\s_1$:
\begin{equation}
0\longrightarrow\s_2\longrightarrow\s\longrightarrow\s_1\longrightarrow0.
 \end{equation}
Note that $\s_1$ is a minimal parabolic subalgebra of $\s\u(1,1)$ while 
$\s_2$ is a minimal parabolic subalgebra of $\s\u(1,n)$. In particular, 
the Lie algebra $\s$ is exact symplectic w.r.t. the element 
$\xi:=\xi_1\oplus\xi_2$ of $\s^\star$ with $\xi_i\in\s_i^\star\quad(i=1,2)$ 
defined as $\xi_1:=N_0^*$ and $\xi_2:=N_3^*$.

One therefore obtains a UDF for proper actions of $\CA\CN$ by direct 
application of the above extension lemma \ref{EXT}. 

\section{Isospectral deformations of anti de Sitter black holes}

\subsection{Anti de Sitter black holes}
Anti de Sitter (AdS) black holes have been introduced by Ba\~nados, Teitelbaum, 
Zannelli and Henneaux \cite{BTZ, BHTZ} as connected locally AdS space-times 
$\CM$ (possibly with boundary and corners) admitting a {\sl singular} causal 
structure in the following sense: 

\begin{condition}
 There exists a closed subset $\CS$ in $\CM$ 
called the \emph{singularity} such that the subset $\CM_{\mbox{bh}}$  
constituted by all the points $x$ such that every light like geodesic issued 
from $x$ ends in $\CS$ within a finite time is a proper open subset of 
$\CM_{\mbox{phys}}:=\CM\backslash\CS$.
\end{condition}

Originally such solutions were constructed in space-time dimension 3, but they   
exist in arbitrary dimension $n \geq 3$ (see \cite{BDRS,ClDe}).                      
More precisely the structure may be described as follows. Take 
$\mathcal{G}:=SO(2,n-1)$ (the AdS group), fix a Cartan involution $\theta$ and 
a $\theta$-commuting involutive automorphism $\sigma$ of $\mathcal{G}$ such that 
the subgroup $\mathcal{H}$ of $\mathcal{G}$ of the elements fixed by $\sigma$ 
is locally isomorphic to $SO(1,n-1)$. The quotient space 
$\mathcal{M}:=\mathcal{G}/\mathcal{H}$ is an $n$-dimensional Lorentzian     
symmetric space, the {\sl anti de Sitter space-time}.          
It is a solution of the Einstein equations without source.                
Let $\g$ denote the Lie algebra of $\mathcal{G}$ and denote by $\g=\h\oplus\q$ 
the $\pm 1$-eigenspace decomposition with respect to the differential at $e$ 
of $\sigma$ that we denote again by $\sigma$. Denote by $\g=\k\oplus\p$ the 
Cartan decomposition induced by $\theta$, consider a $\sigma$-stable    
maximally Abelian subalgebra $\a$ in $\p$ and choose accordingly a positive 
system of roots. Denote by $\n$ the corresponding nilpotent subalgebra. 
Set  $\overline{\n}:=\theta(\n)$, $\r:=\a\oplus\n$ and 
$\overline{\r}:=\a\oplus\overline{\n}$. Finally denote by $\CR:=\CA\CN$ and 
$\overline{\CR}:=A\overline{N}$ the corresponding analytic subgroups 
of $\mathcal{G}$. One then has 
\begin{proposition}\cite{BDRS,ClDe} 
The groups $\CR$ and $\overline{\CR}$ admit open orbits and finitely many 
closed orbits in the AdS space $\mathcal{M}$.  Prescribing as singular the 
union of all closed orbits (of $\CR$ and $\overline{\CR}$) defines a structure 
of causal black hole on an open subset $\CM_{\mbox{phys}}$ in $\mathcal{M}$ 
(in the sense of the above condition {\bf (BH)}).  In particular, every open 
orbit $\CM_{o}$ of $\CR$ in $\mathcal{M}$ containing $\CM_{\mbox{phys}}$ is 
itself endowed with a black hole structure. 
\end{proposition}

Recall that if $J$ denotes an element of $Z(\CK)$ whose associated conjugation 
coincides with the Cartan involution $\theta$ then the $\CR$-orbit $\CM_o$ in   
$\mathcal{G}/\mathcal{H}$ of an element $u\mathcal{H}$ with $u^2=J$ is open 
and contains $\CM_{phys}$, see \cite{ClDe}. Remark that the extension 
lemma \ref{EXT} yields an oscillatory integral UDF for proper actions 
of $\CR$. But here the situation is simplified by the following         
observation --- for convenience of the presentation we write it below   
for $n=4$ but the results are valid for any $n\geq 3$.                  

\begin{proposition}\label{GT}
The $\CR$-homogeneous space $\CM_o$ admits a unique structure of globally 
group type symplectic symmetric space. The latter is isomorphic to 
$(\CR_0,\omega,s)$ described in section 2.
\end{proposition}

For the proof, recall first (see \cite{ClDe}) that the solvable part of the
Iwasawa decomposition of $\so(2,3)$ may be realized with as nilpotent part
$\n=\{W,V,M,L\}$ and Abelian $\a=\{J_1,J_2\}$  
with the commutator table
$[V,W]=M$,	$[V,L]=2W$,
$[J_1,W]=W$,	$[J_2,V]=V$,
$[J_1,L]=L$,	$[J_2,L]=-L$,
$[J_1,M]=M$,	$[J_2,M]=M$.
Notice that $W,J_{1}\in\h$, and $J_{2}\in\q$. This decomposition is related 
to the one given in \eqref{EqTableSP} by 
 \begin{align*}
N_0&=L&N_2&=2W&H_1&=J_1+J_2\\
N_1&=V&N_3&=M&H_2&=J_1-J_2.
\end{align*}
We choose to study the orbit of the element $\vartheta=u\mathcal{H}$ with
\[ 
  u=\begin{pmatrix}
0&1\\
-1&0\\
&&I_{3\times 3}
\end{pmatrix}.
\]
We denote by $\CR_{\theta}$ its stabilizer group in $\CR$, by $\r_{\theta}$ the
Lie algebra of $\CR_{\theta}$; $\r'$ is the subalgebra of $\r$ generated by 
elements of $\r$ minus the generator of $\r_{\theta}$ and $\CR'$ is the 
analytic subgroup of $\CR$ whose algebra is $\r'$.
\begin{lemma}
The action of $\CR'$ on $\mU$ is simply transitive, i.e. 
$\CR'u\mathcal{H}=\CR u\mathcal{H}$.
\end{lemma}

\begin{proof}
The first step is to prove that $\CR_{\theta}$ is connected and $\mU$ 
simply connected in order to prevent any double covering problem. 
The stabilizer of $ u\mathcal{H}$ is
\begin{equation}   
\CR_{\theta}=\{ r\in \CR\mid r\cdot u\mathcal{H}= u\mathcal{H} \}
	=\{ r\in \CR\mid C_{u^{-1}}(r)\in \mathcal{H} \}.     
\end{equation}
Since $\CR$ is an exponential group, we have 
  $\r_{\theta}=\{ X\in\CR\mid\Ad(u^{-1})X\in\h \}$ with
$\CR_{\theta}=\exp\r_{\theta}$.
The set $\r_{	\theta}$ being connected, $\CR_{\theta}$ is connected too.
A long exact sequence argument using the fibration $\CR_{\theta}\to \CR\to\mU$ 
shows that $H^0(\CR_{\theta})\simeq H^1(\mU)$, which proves that $\mU$ is 
simply connected. 

As an algebra, $\r$ is a split extension $\r=\r_{\theta}\oplus_{\ad}\r'$. 
Hence, as group, $\CR=\CR_{\theta} \CR'$, or equivalently 
$\CR=\CR'\CR_{\theta}$. This proves that the action is transitive. 
The action is even simply transitive because $\mU$ is simply connected.
\end{proof}

Let us now find the algebra $\r_{\theta}$. The Cartan involution 
$X\mapsto -X^{t}$ is implemented as $C_{J}$ with
\[ 
  J=\begin{pmatrix}
-I_{2\times 2}\\
&I_{3\times 3}
\end{pmatrix}.
\]
Using the relations $u^{2}=J$ and $\sigma(u)=u^{-1}$, one sees that 
$C_{u^{-1}}(r)\in\CH$ if and only if 
$\sigma\Big( C_{u^{-1}}r \Big)=C_{u^{-1}}r$. 
This condition is equivalent to                  
$\theta\sigma (r)=r$.                            
The involution $\sigma$ splits $\a$ into two parts: 
$\a=\a^{+}\oplus\a^{-}$ with $J_1\in\a^{+}=\a\cap\h$ and 
$J_2\in\a^{-}\cap\q$. Let $\beta_{1},\beta_{2}\in\a^*$ be the dual basis; 
we have $W\in\g_{\beta_{1}}$, $V\in\g_{\beta_{2}}$, 
$L\in\g_{\beta_{1}-\beta_{2}}$, $M\in\g_{\beta_{1}+\beta_{2}}$, and, 
in terms of positive roots, the space $\n$ is given by
$W\in\g_{\alpha+\beta}$, $V\in\g_{\beta}$,
$L\in\g_{\alpha}$, $M\in\g_{\alpha+2\beta}$.
We are now able to compute the vectors $X\in\r$ such that $\sigma\theta(X)=X$. 
Let us take $X\in\r=\a\oplus\n$ and apply $\sigma\theta$:
\begin{equation}
\begin{split}
X&=X_{J_1}+X_{J_2}+X_{\alpha}+X_{\beta}+X_{\alpha+\beta}+X_{\alpha+2\beta},\\
\sigma\theta X&=-X_{J_1}+X_{J_2}+Z_{-(\alpha+2\beta)}+
Z_{\beta}+Z_{-(\alpha+\beta)}+Z_{-\alpha}
\end{split}
\end{equation}
where $X_{\varphi}$ and $Z_{\varphi}$ denote elements of $\g_{\varphi}$. 
It is directly apparent that $X=J_2$ belongs to $\r_{\theta}$. The only other 
component common to $X$ and $\sigma\theta X$ is in $\g_{\beta}$, but   
it is \emph{a priori} not clear that $X_{\beta}=Z_{\beta}$. The dimension 
of $\mU$ is $4$ and that of $\CR$ is $6$, hence $\CR_{\theta}$ is at least  
$2$-dimensional; it is generated by $J_2$ and $\g_{\beta}=\R V$, i.e.    
$\r_{\theta}=\Span\{ J_2,V \}$.                       
This proves that the orbit of $ u\mathcal{H}$ is open.          

The fact that $\CR'$ acts freely on $\mU=\CR/\CR_{\theta}$ proves that $\mU$ 
is locally of group type and since, by definition, $\mU$ is only one orbit 
of $\CR$, the space $\mU$ is globally of group type. From now on, 
$\CM_o=\CR/\CR_{\theta}$ will be identified with $\mU$ as homogeneous space, 
so what we have to find is a group $\tilde{\CR}$ which
\begin{itemize} 
\item  acts transitively on $\mU$, i.e. $\tilde{\CR} u\mathcal{H}=\CR u\mathcal{H}$,
\item admits a symplectic structure.
\end{itemize}
It is immediate to see that the algebra $\r'$ fails to fulfil the symplectic 
condition. The algebra $\tilde{\r}=\Span\{ A,B,C,D \}$ of a group which 
fulfils the first condition must at least act transitively on a small    
neighborhood of $ u\mathcal{H}$ and thus be of the form
\begin{subequations}   \label{EqAlgGEnennsy}
\begin{align} 
 A&=J_{1}+aJ_{2}+a'V\\
 B&=W+bJ_{2}+b'V\\
 C&=M+cJ_{2}+c'V\\
 D&=L+dJ_{2}+d'V.
\end{align}
\end{subequations}
The problem is now to fix the parameters $a,a',b,b',c,c',d,d'$ in such a way 
that $\Span\{ A,B,C,D \}$ is a Lie algebra (i.e. it is closed under the Lie 
bracket) which admits a symplectic structure and whose group acts transitively 
on $\mU$.  We will begin by proving that the surjectivity condition imposes   
$b=c=d=0$. Then the remaining conditions for $\tilde \r$ to be an algebra 
are easy to solve  by hand.

First, remark that $A$ acts on the algebra $\Span\{ B,C,D \}$ because $J_1$ 
does not appears in $[\r,\r]$. We can write 
$\tilde{\r}=\R A\oplus_{\ad}\Span\{ B,C,D \}$                      
and therefore a general element of the group $\tilde{\CR}$ reads 
$\tilde r (\alpha,\beta,\gamma,\delta)= e^{\alpha A} e^{\beta B+\gamma C+\delta D}$ 
because a subalgebra of a solvable exponential Lie algebra is solvable exponential. 
Our strategy will be to split this expression in order to get a product 
$\CS\CR'$ (which is equivalent to a product $\CR'\CS$). As Lie algebras, 
$\Span\{ B,C,D \}\subset\R J_2\oplus_{\ad}\{ W,M,L,V \}$. Hence there exist 
functions $w$, $m$, $l$, $v$ and $x$ of $(\alpha,\beta,\gamma,\delta)$ such that
\begin{equation} \label{EqGeneRi}
 e^{\beta B+\gamma C+\delta D}= e^{xJ_2} e^{wW+mM+lL+vV}.
\end{equation}
We are now going to determine $l(\alpha,\beta,\gamma,\delta)$ and study the 
conditions needed in order for $l$ to be surjective on $\R$. Since $J_2$ does not 
appear in any commutator, the Campbell-Baker-Hausdorff formula yields 
$x=\beta b+\gamma c+\delta d$. From the fact that $[J_2,L]=-L$, we see that 
the coefficient of $L$ in the left hand side of \eqref{EqGeneRi} is 
$-l(1- e^{-x})/x$. The $V$-component in the exponential can also get out 
without changing the coefficient of $L$. We are left with     
$\tilde r(\alpha,\beta,\gamma,\delta)= e^{\alpha A} e^{xJ_2} e^{yV} e^{w'W+m'M+lL}$ 
where $w'$ and $m'$ are complicated functions of $(\beta,\gamma,\delta)$ 
and $l$ is given by
\begin{equation}
l(\beta,\gamma,\delta)=\frac{ -\delta (\beta b+\gamma c+\delta d) }
{ 1- e^{-\beta b-\gamma c-\delta d} },
\end{equation}
which is not surjective except when $b=c=d=0$. Taking the inverse a general 
element of $\tilde\CR u\mathcal{H}$ reads 
$\big[  e^{ -wW-mM-lM} e^{j_1J_1}u \big]$, where the range of $l$ is not 
the whole $\R$.  Since the action of $R'$ is \emph{simply} transitive,     
$\tilde\CR$ is not surjective on $\CR u\mathcal{H}$.

When $b=c=d=0$, the conditions for \eqref{EqAlgGEnennsy} to be an algebra are 
easy to solve, leaving only two \emph{a priori} possible two-parameter 
families of algebras:
\paragraph{Algebra 1}  
\begin{subequations}
 \begin{align*}
A&=J_{1}+\frac{ 1 }{2}J_{2}+sV	&[A,B]&=B+sC\\
B&=W				&[A,C]&=\frac{ 3 }{2}C\\
C&=M				&[A,D]&=2sB+\frac{ 1 }{2}D\\
D&=L+rV				&[B,D]&=-rC.
\end{align*}
\end{subequations}
with $r\neq 0$. The general symplectic form on that algebra is given by  
\begin{equation}   
\omega_{1}=\begin{pmatrix}
0	&-\alpha		&-\beta	&-\gamma\\
\alpha	&0			&0	&\frac{ 2\beta r }{ 3 }\\
\beta	&0			&0	&0\\
\gamma	&-\frac{ 2\beta r }{ 3 }	&0	&0
\end{pmatrix},
\end{equation}                                               
Since $\det\omega=\left( \frac{ 2\beta r }{ 3 } \right)^{2}$    
we must have $\beta\neq 0$, $r\neq 0$. That algebra will be denoted by 
$\r_{1}$. The analytic subgroup of $\CR$ whose Lie algebra is $\r_1$ is 
denoted by $\CR_1$.

\paragraph{Algebra 2}
\begin{align*}
A&=J_{1}+rJ_{2}+sV		&[A,B]&=B+sC\\
B&=W				&[A,C]&=(r+1)C\\
C&=M				&[A,D]&=2sB+(1-r)D\\
D&=L.
\end{align*}
There is no way to get a non-degenerate symplectic form on that algebra. 
\begin{remark}
One can eliminate the two parameters in algebra $\r_{1}$ by the isomorphism  
\begin{equation}
\phi=
\begin{pmatrix}
1&0&0&0\\
0&1&0&4s\\
0&2sr&1/r&4s^{2}/r\\
0&0&0&1
\end{pmatrix}
\end{equation}
which fixes $s=0$ and $r=1$ and transforms $\r_{1}$ into the algebra defined by 
  $[A',B']=B'$,  $[A',C']=\frac{ 3 }{2}C'$,
$[A',D']=\frac{ 1 }{2}D'$, $[B',D']=-C'$.
\end{remark}
It is now easy to prove that 
\begin{proposition}   \label{PropCRunXXX}
The group $\CR_1$ of algebra $\r_1$ acts transitively on $\mU$, 
i.e. $\CR u\mathcal{H}=\CR_1 u\mathcal{H}$.
\end{proposition}

\begin{proof}
The algebra $\r_1$ can be written 
  $\r_1=\R A\oplus_{\ad}\R D\oplus_{\ad}\Span\{ B,C \}$, a split extension,
hence a general element reads
$r_1(\alpha,\beta,\gamma,\delta)= e^{\alpha A} e^{\delta D} e^{\beta W+\gamma M}$.
One can use Campbell-Baker-Hausdorff formula to split it into a factor in 
$\CR_{\theta}$ and one in $\CR'$ (where $f$ and $g$ are some auxiliary functions):
\begin{equation}
r_1(\alpha,\beta,\gamma,\delta)=\underbrace{e^{\alpha s V+
\frac{ \alpha }{2}J_2} e^{\delta rV}}_{\in\CR_{\theta}} 
\underbrace{e^{\alpha J_1} e^{\big( f(\delta)+\beta \big)W\big( g(\delta)+
\gamma \big)M}}_{\text{surjective on $\CR'$}}
\end{equation}

\end{proof}

The conclusion is that $\CR_1$ is the group $\tilde{\CR}$ that we were 
searching for. To summarize, the structure is as follows. 

\begin{enumerate}
\item 
The AdS space is decomposed into a family of cells: the orbits of a symplectic 
solvable Lie group $\tilde{\CR}$ as in Proposition~\ref{GT} above.  
Note that these cells may be viewed as the symplectic leaves of the Poisson 
generalized foliation associated with the left-invariant symplectic 
structure on $\tilde{\CR}$.
\item 
The open $\tilde{\CR}$-orbit $\CM_o$, endowed with a black hole structure, 
identifies with the group manifold $\tilde{\CR}$.
\end{enumerate}

\subsection{Deformation triples for $\CM_o$}

\subsubsection{Left-invariant Hilbert function algebras on $\CR_0$}

In this section, we present a modified version of the oscillatory integral 
product (\ref{PRODUCT}) leading to a left-invariant associative algebra 
structure on the space of square integrable functions on $\CR_0$.

\begin{theorem} \label{thmL2} Let $u$ and $v$ be smooth compactly supported 
functions on $\CR_0$. Define the following three-point functions:
 \begin{equation}
\begin{split}
S:=& S_V\big(\cosh(a_1-a_2)x_0, \cosh(a_2-a_0)x_1, \cosh(a_0-a_1)x_2\big)\\
&-\cyclic_{0,1,2}\sinh\big(2(a_0-a_1)\big)z_2\;;
\end{split}
\end{equation}
and 
\[ 
\begin{split}
A:= \Big[&\cosh\big(2(a_1-a_2)\big)\cosh(2(a_2- a_0))\cosh(2(a_0-a_1))\\
& \big[\cosh(a_1-a_2)\cosh(a_2- a_0)\cosh(a_0-a_1)\big]^{\dim\CR_0-2}\Big]^{\frac{1}{2}}.
\end{split}
\]
Then the formula                                                       
\begin{equation}\label{HILB}
u\star^{(2)}_\theta v\;:=\;\frac{1}{\theta^{\dim {\mathcal R}_0}}
\int_{\CR_0\times\CR_0}A\;e^{\frac{2i}{\theta}S}u\otimes v
\end{equation}
extends to $L^2(\CR_0)$ as a left-invariant associative Hilbert algebra 
structure. In particular, one has the strong
closedness\footnote{The notion of strongly closed star product was 
introduced in  \cite{CFS} in the formal context.} property:
\begin{equation*}
\int u\star^{(2)}_\theta v\;=\;\int uv\;.      
\end{equation*}
\end{theorem}
\begin{proof}
The oscillatory integral product (\ref{PRODUCT}) may be obtained by 
intertwining the Weyl product on the Schwartz space $\CS$ (in the Darboux global 
coordinates \eqref{DARBOUX}) by the following integral operator \cite{BiMas}:
\begin{equation*}
\tau:=\CF^{-1}\circ(\phi_\theta^{-1})^\star\circ\CF\;,
\end{equation*}
$\CF$ being the partial Fourier transform with respect to        %
the central variable $z$:
\begin{equation*}
\CF(u)(a,x,\xi):=\int e^{-i\xi z}u(a,x,z){\rm d}z\;;
\end{equation*}
and $\phi_\theta$ the one parameter family of diffeomorphism(s):      
\begin{equation*}
\phi_\theta(a,x,\xi)=(a,\frac{1}{\cosh(\frac{\theta}{2}\xi)}x,
\frac{1}{\theta}\sinh(\theta\xi)).
\end{equation*}
Set $\mJ:=|(\phi^{-1})^\star\mbox{Jac}_\phi|^{-\frac{1}{2}}$ and observe that 
for all $u\in C^\infty\cap L^2$, the function $\mJ\,(\phi^{-1})^\star u$ 
belongs to $L^2$. Indeed, one has 
\begin{equation*}
\int|\mJ\,(\phi^{-1})^\star u|^2=
\int|\phi^\star \mJ|^2\,|\mbox{Jac}_\phi|\,|u|^2=\int|u|^2\;.
\end{equation*}
Therefore, a standard density argument yields the following isometry:
\begin{equation*}
T_\theta:L^2(\CR_0)\longrightarrow L^2(\CR_0):
u\mapsto\CF^{-1}\circ m_\mJ\circ(\phi^{-1})^\star\circ\CF(u)\;,
\end{equation*}
where  $m_{\mJ}$ denotes the multiplication by $\mJ$.  
Observing that $T_\theta=\CF^{-1}\circ m_{\mJ}\circ\CF\circ\tau$, 
one has $\star^{(2)}_\theta=\CF^{-1}\circ m_{\mJ}\circ\CF(\star_\theta)$. 
A straightforward computation (similar to the one in \cite{Bi})     
then yields the announced formula.
\end{proof}

\begin{remark}
Let us point out two facts with respect to the above formulas:
\begin{enumerate}
\item Note the cyclic symmetry of the oscillating three-point kernel 
$A\;e^{\frac{2i}{\theta}\,S}$.
\item The above oscillating integral formula gives rise to a 
strongly closed, symmetry invariant, formal star product on the symplectic 
symmetric space $(\CR_0,\omega,s)$.
\end{enumerate}
\end{remark}

\begin{proposition}
The space $L^2(\CR_0)^{\infty}$ of smooth vectors in $L^2(\CR_0)$ of the left 
regular representation closes as a subalgebra of $(L^2(\CR_0),\star^{(2)}_\theta)$.
\end{proposition}

\begin{proof}
First, observe that the space of smooth vectors may be described  
as the intersection of the spaces $\{V_n\}$ where  $V_{n+1}:=(V_n)_1$,   
with $V_0:=L^2(\CR_0)$ and $(V_n)_1$ is defined as the space of elements $a$ of 
$V_n$ such that, for all $X\in\r_0$, $X.a$ exists as an element of $V_n$ 
(we endow it with the projective limit topology).                      

Let thus $a,b\in V_1$. Then, $(X.a)\star b+a\star(X.b)$ belongs to $V_0$. 
Observing that $\CD\subset V_1$ and approximating $a$ and $b$ by sequences 
$\{a_n\}$ and $\{b_n\}$ in $\CD$, one gets (by continuity of $\star$): 
$(X.a)\star b+a\star(X.b)=\lim(X.a_n\star b_n+a_n\star(X.b_n))=
\lim X.(a_n\star b_n)=X.(a\star b)$. Hence $a\star b$ belongs to $V_1$. 
One then proceeds by induction.                                       
\end{proof}

 \subsubsection{Twisted $L^2$-spinors and deformations of the Dirac operator}

We now follow in our four dimensional setting the deformation scheme 
presented in \cite{BDRS} in the three-dimensional BTZ context.

At the level of the (topologically trivial) open $\CR$-orbit,  the spin 
structure over $\CM_o$ and the associated spinor $\C^2$-bundle 
-- restriction of the spinor bundle on AdS$_n$ -- are trivial.   
The space of (smooth) spinor fields may then be viewed as
${\bf S}\;:=\;C^\infty(\CR_0,\C^2)$, on which the (isometry) group 
$\CR_0$ acts on via the left regular representation.
In this setting, the restriction to $\CM_o$ of the Dirac operator $D$ 
on AdS$_n$ may be written as                               
$D\;=\;\sum_i\,\gamma^i\,(\tilde{X_i}\,+\,\Gamma_i)$, where
\begin{enumerate}
\item[$\bullet$] $\{X_i\}$ denotes an orthonormal basis of 
$\r_0=T_\vartheta(\CM_o)$ (w.r.t. the adS-metric at the base point $\vartheta$
of $\CM_o$);
\item[$\bullet$] for $X\in\r_0$, $\tilde{X}$ denotes the associated 
left-invariant vector field on $\CR_0$;
\item[$\bullet$] $\gamma^i$ and $\Gamma_i$ denote respectively the Dirac 
$\gamma$-endomorphism and the spin-connection element associated with $X_i$.
\end{enumerate}
In that expression the elements $\gamma^i$'s and $\Gamma_i$'s are constant. 
However, already at the formal level, a left-invariant vector field $\tilde{X}$ 
as infinitesimal generator of the right regular representation does not 
in general act on the deformed algebra.    
In order to cure this problem, we twist the spinor module in the following way.

\begin{definition}
Let ${\rm d}^rg$ be a right-invariant Haar measure on $\CR_0$ and consider the 
associated space of square integrable functions $L^2_{\mathrm{right}}(\CR_0)$. Set
\begin{equation*}
\CH\;:=\;L^2_{\mathrm{right}}(\CR_0)\otimes\C^2\;;
\end{equation*}
and denote by $\CH^\infty$ the space of smooth vectors in $\CH$ of the 
natural right representation of $\CR_0$ on $\CH$. Then intertwining 
$\star^{(2)}_\theta$ by the inverse mapping 
$$\iota:L^2_{\mathrm{right}}(\CR_0)\to L^2(\CR_0):\iota(u)(g):=u(g^{-1})$$ 
yields a right invariant noncommutative 
$L^2_{\mathrm{right}}(\CR_0)$-bi-module structure (respectively a         
$(L^2_{\mathrm{right}}(\CR_0))^\infty$-bi-module structure) on $\CH$ 
(resp. $\CH^\infty$).  The latter will be denoted $\CH_\theta$ 
(resp. $\CH^\infty_\theta$).
\end{definition}

We summarize the main results of this paper in the following:

\begin{theorem}
The Dirac operator $D$ acts in $\CH^\infty_\theta$ as a derivation of the 
noncommutative bi-module structure. In particular, for all 
$a\in (L^2_{\mathrm{right}}(\CR_0))^\infty$, the commutator $[D,a]$ extends 
to $\CH$ as a bounded operator. In other words, the triple 
$(L^2_{\mathrm{right}}(\CR_0))^\infty, \CH^\infty_\theta, D)$ induces 
on $\CM_o$ a pseudo-Riemannian deformation triple.                  
\end{theorem}

\section{Conclusions, remarks and further perspectives}

To the AdS space we associated a symplectic symmetric space $(M,\omega,s)$.
That association is \textit{natural} by virtue of the uniqueness property
mentioned in Proposition~\ref{GT}. The data of any invariant (formal or not)
deformation quantization on $(M,\omega,s)$ yields then canonically a UDF 
for the actions of a non-Abelian solvable Lie group. Using it we defined the 
noncommutative Lorentzian spectral triple $({\mathcal A}^\infty,\CH,D)$ where 
${\mathcal A}^\infty:=(L^2_{\mathrm{right}}(\CR_0))^\infty$ is a noncommutative
Fr\'echet algebra modelled on the space of smooth vectors of the regular 
representation on square integrable functions on the group $\CR_0$. 
The underlying commutative limit is endowed with a causal black hole structure 
encoded in the $\CR_0$-group action. 
A first question that this construction raises is that of defining within 
the present Lorentzian context the notion of causality at the 
operator algebraic level.

Another direction of research is to analyze the relation between the present 
geometrical situation and the corresponding one within the quantum group 
context. Indeed, our universal deformation formulas can be used at the 
algebraic level to produce nonstandard quantum groups $SO(2,n-1)_q$ via 
Drinfeld twists. An interesting challenge would then be to study the 
behaviour of the representation theory under the deformation process.

More generally the somewhat elliptic sentence with which we started the paper
may now be better understood if we remark that the physical motivation 
section and the quantum group framework suggest to study a number of 
questions related to (noncommutative) singleton physics, in particular:
\begin{enumerate}

\item Since \cite{FHT93,Sta98} we know that for $q$ even root of unity, 
there are unitary irreducible finite dimensional representations of the
Anti de Sitter groups. Interestingly (cf. \cite{FGR99} p.122) the    
``fuzzy 3-sphere" is related to the Wess-Zumino-Witten models 
and is conjectured to be related to the non-commutative geometry of 
the quantum group $U_q({\rm sl}_2)$ for $q=e^{2\pi i/(k+2)}, k>0$, 
a root of unity. 

\item The last remark suggests to look more closely at the phenomenon
of dimensional reduction which appears in a variety of related problems.
In this paper we considered only $n \geq 3$. The reason is that for
$n=2$ the context is in part different: the conformal group of 1+1
dimensional space-time is infinite dimensional, and there are no black
holes \cite{ClDe}. But many considerations remain true, and furthermore
many group-theoretical properties find their origin at the 1+1 dimensional
level, e.g. the uniqueness of the extension to conformal group \cite{AFFS}.
Another exemple of dimensional reduction is the fact that the massless
UIR of the 2+1 dimensional Poincar\'e group Di and Rac satisfy
Di $\oplus$ Rac = D(HO)$\oplus$D(HO) where D(HO) is the representation
D(1/4)$\oplus$D(3/4) of the metaplectic group (double covering of 
$SL(2,\mathbb{R})$) which is the symmetry of the harmonic oscillator in the 
deformation quantization approach (see e.g. Section (2.2.4) in \cite{DS02}).     

\item What do the degenerate representations Di and Rac become under 
deformation? Furthermore there may appear, for our nonstandard quantum
group $SO(2,n)_q$, new representations that have no equivalent at the
undeformed level (e.g. in a way similar to the supercuspidal representations
in the $p$-adic context). These may have interesting physical interpretations.

\item We have seen that for $q$ even root of unity $SO(2,n)_q$ has some
properties of a compact Lie group. Our cosmological Ansatz suggests that
the $q$AdS black holes are ``small." It is therefore natural to try and
find a kind of generalized trace that permits to give a finite volume for
$q$AdS. Note that, in contradistinction with infinite dimensional Hilbert
spaces, the notions of boundedness and compactness are the same for closed
sets in Montel spaces, and that our context is in fact more Fr\'echet nuclear
than Hilbertian. This raises the more general question to define in an
appropriate manner the notion of ``$q$-compactness" (or ``$q$-boundedness")
for noncommutative manifolds. 

\item Possibly in relation with the preceding question, one should perhaps
consider deformation triples in which the Hilbert space is replaced 
by a suitable locally convex topological vector space (TVS), on which
$D$ could be continuous.  

\item The latter should yield a natural framework for implementing quantum
symmetries in deformation triples, since Fr\'echet nuclear spaces and their
duals are at the basis of the topological quantum groups (and their duals)
introduced in the 90's, especially in the semi-simple case with preferred
deformations (see the review \cite{BGGS}). We would thus in fact have
quadruples $(\mathcal{A}, \mathcal{E}, D, \mathcal{G})$ where $\mathcal{A}$ is
some topological algebra, $\mathcal{E}$ an appropriate TVS, $D$ some 
(bounded on $\mathcal{E}$) ``Dirac" operator and $\mathcal{G}$ some symmetry. 
[Being in a Lorentzian noncompact framework, we did not address here 
questions such as the resolvent of $D$ when $\mathcal{E}$ is a Hilbert      
space, which we did not need at this stage; eventually one may however
have to deal with the reasons that motivated the additional requirements on
triples in the Riemannian compact context; note that here the restriction to
an open orbit was needed in order to have bounded commutators $[D,a]$ in the
Hilbertian context, but a good choice of $\mathcal{E}$ could lift the
restriction.] That framework should be naturally extendible to the
supersymmetric context, which is the one considered in \cite{FGR99}
with modified spectral triples and is natural also for the problems
considered here since e.g. Di$\oplus$Rac and D(HO) are UIR of the
corresponding supersymmetries.

\item If we want to incorporate ``everything," the (external) symmetry
$\mathcal{G}$ should be the Poincar\'e group $SO(1,3)\cdot \mathbb{R}^4$ in the
ambient Minkowski space (possibly modified by the presence of matter)
and $SO(2,3)_q$ in the $q$AdS$_4$ black holes, or possibly some
supersymmetric extension. The unified (external) symmetry could therefore
be something like a groupoid. The latter should be combined in a subtle way
(as hinted e.g. in \cite{St07}) with the ``internal" symmetry associated with
the various generations, colors and flavors of (composite) ``elementary"
particles in a generalized Standard Model, possibly in a noncommutative
geometry framework analogous to what is done in \cite{Co06,CCM,Ba06}.     
There would of course remain the formidable task to develop quantized 
field theories on that background, incorporating composite QED for photons 
on AdS as in \cite{FF88} and some analog construction for the electroweak model
(touched in part in \cite{Fr00}) and for QCD, possibly making use of some 
formalism coming from string theory. 

\item The Gelfand isomorphism theorem permits to realize commutative
involutive algebras as algebras of functions on their ``spectrum."
Finding a noncommutative analog of it has certainly been in the back of the
mind of many, since quite some time (see e.g. \cite{St05}).
We now have theories and many examples of deformed algebras, quantum groups
and noncommutative manifolds. The above mentioned quadruples could
provide a better understanding of that situation. 
 
\end{enumerate}

\bibliographystyle{amsalpha}

\end{document}